\documentstyle[amscd,amssymb,verbatim,12pt,righttag]{amsart}
\numberwithin{equation}{section}

\theoremstyle{plain}
\newtheorem{thm}{Theorem}[section]
\newtheorem{lemma}[thm]{Lemma}
\newtheorem{pro}[thm]{Proposition}
\newtheorem{cor}[thm]{Corollary}
\newtheorem{ex}[thm]{Example}
\newtheorem{de}[thm]{Definition}
\newtheorem{re}[thm]{Remark}

\begin{document}
%Topmatter
\title[Ruelle Operator Theorem]%
{Ruelle Operator Theorem for Nonexpansive systems}
 \author[YunPing Jiang and Yuan-Ling Ye ]{}

\thanks{** The corresponding author.\\
\indent * Supported by PSC-CUNY awards and the CUNY Collaborative
Research Incentive Program
and the Bai Ren Ji Hua from Academia Sinica.\\
\indent ** Supported by the Project-sponsored by SRF for ROCS, SEM.}

%End Topmatter

\maketitle

\begin{center}
\author{Yunping Jiang* and Yuan-Ling Ye**}
\end{center}

\bigskip

\begin{center}
Abstract
\end{center}

The Ruelle operator theorem has been studied extensively both in
dynamical systems and iterated function systems. In this
paper we study the Ruelle operator theorem for nonexpansive systems. Our theorems give some sufficient conditions
for the Ruelle operator theorem to be held for a nonexpansive system.

%%%%%%%%%%%%%%%%%%

%%%%%%%%%%%%%%%%%%
\bigskip
\bigskip

\section{\bf Introduction }

Ruelle introduced a convergence theorem to study the equilibrium
state of an infinite one-dimensional lattice gas in his famous
paper~\cite{Rue}. Bowen~\cite{Bow} further set up the theorem as
the convergence of powers of a Ruelle operator on the space of
continuous functions on a symbolic space. More precisely, let
$$
\Sigma = \{1, \cdots N\}^{\Bbb N}=\{ \omega =i_{0}i_{1}\cdots
i_{n-1} \cdots \;|\; i_{n-1}\in \{ 1, \cdots, N\},\; n=1,
2,\cdots\}
$$
be the one-sided symbolic space and
$$
\sigma: \omega =i_{0}i_{1}\cdots i_{n-1} \cdots \to \sigma
(\omega) =i_{1}\cdots i_{n-1} \cdots
$$
be the left shift of $\Sigma$. Then $(\Sigma, \sigma)$ is called a
symbolic system. Let $\phi$ be a H\"older continuous function on
$\Sigma$ (a potential). Let $C(\Sigma)$ be the space of all
continuous functions on $\Sigma$. The {\it Ruelle operator} is
defined as
\begin{equation}
 {\cal T}f (x) = \sum_{y \in \sigma^{-1}(x)} e^{\phi (y)} f(y), \qquad  f \in
C(\Sigma).
\end{equation}
It is a positive operator, that is, ${\cal T}f >0$ whenever $f>0$.

Let $\varrho$ be the spectral radius of the operator
$$
{\cal T}: C(\Sigma)\to C(\Sigma).
$$
It is known that $\varrho$ is the unique positive simple maximal
eigenvalue of ${\cal T}$ acting on the space of all H\"older
continuous functions on $\Sigma$ (see, for example,~\cite{Jia}).
It was then proved that ${\cal T}$ has a unique positive
eigenfunction $h \in C(\Sigma)$ and a unique probability
eigenmeasure $\mu \in C^{\ast} (\Sigma)$ corresponding to the
eigenvalue $\varrho>0$ (see, for example,~\cite{Bow}). And
moreover, for any $f \in C(\Sigma)$, $\varrho^{-n}{\cal T}^n(f)$
converges uniformly to a constant multiple of $h$. This is called
the Ruelle operator theorem. In this theorem, $\sigma: \Sigma\to
\Sigma$ is an expanding dynamical system. More general results
about the Ruelle operator theorem for expanding dynamical systems
and contractive iterated function systems (IFS) have been also
obtained. We give a partial list in the
literature~\cite{Fan,FaJi1,FaJi2,FaLa,Wa1,Wa2}.

Recently a parabolic system has drawn a great attention to people
who are interested in the Ruelle operator theorem (refer to~\cite
{BJL,LSV,MaUr,PrSl,Urb,Ye1,Ye2,You,Yur}). However, in this case,
it is known that the bounded eigenfunction of the spectral radius
$\varrho$ of $T$ may not exist~\cite{LaYo}, and even if the
eigenfunction exists, $\varrho$ may not be an isolated point of
the spectrum~\cite{BDEG}. So far the results known are far from
satisfactory. And a study of such a system remains a challenge
problem. Lau and Ye studied the Ruelle operator theorem for a
nonexpansive system in a recent paper~\cite{LaYe}. In this paper
we continue to study the above mentioned problem for a
nonexpansive system. In the paper~\cite{LaYe}, one requirement is
that one of the iterations of the IFS must be strictly
contractive. It is important to remove this requirement because
many examples of IFS will not satisfy this requirement. In this
paper, we remove this requirement. It is a major improvement.

Our iterated function system (IFS) $\{w_j\}_{j=1}^m$ in this paper
is weakly contractive as defined by
$$
\alpha_{w_j}(t) :=\sup_{|x-y | \leq t}|w_j(x)- w_j(y)| < t,\quad
\forall t > 0, \quad 1\leq j\leq m
$$
or, more generally, nonexpansive as defined by
$$
|w_j(x) -w_j(y)| \leq |x-y|, \quad 1\leq j \leq m.
$$
For the weakly contractive case, the invariant compact set $K$
exists as in the contractive case (Hata~\cite{Hat}). For the
nonexpansive case we can take the smallest compact invariant $K$
(see Proposition 2.1 for the additional assumption). With each
$w_{j}$, we associate a positive continuous function $p_j$
as a weight function (or potential function).
We can set up the Ruelle operator as in (1.2) on the
space $C(K)$ of continuous functions on $K$,
\begin{equation}
T(f)(x) = \sum_{j=1}^m p_{j}(x)f(w_j(x)),  \qquad  f \in C(K).
\end{equation}
Let $\varrho$ still be the spectral radius of the operator
$$
T: C(K)\to C(K).
$$

\begin{de}~\label{Def1} We call $(X, \{w_{j}\}_{j=1}^{m}, \{p_{j}\}_{j=1}^{m})$ a
nonexpansive system, if all maps $w_j$ are nonexpansive and all
potentials $p_{j}(x)$ are Dini continuous on $X$.
\end{de}

The main result in this
paper which we are particularly interested in is that

\bigskip
\begin{thm}[Main Theorem]~\label{mainth}  {\em Let $(X, \{ w_j\}_{j=1}^{m},
\{ p_j\}_{j=1}^{m})$ be a nonexpansive system. Suppose
$$
\sup_{x \in K} \sum_{j=1}^m p_j (x) \sup_{y \not = x}
\frac{|w_j(x)- w_j(y)|}{|x-y|} < \varrho.
$$
Then the Ruelle operator theorem holds for this nonexpansive
system.}
\end{thm}

\bigskip

We will prove a more general result (Theorem~\ref{3.5}) in \S4.
Actually, the above theorem is a special case of this more general
result. The results in this paper extend the results
in~\cite{LaYe}. However, as we pointed out before, it is a
non-trivial generalization: In the paper ~\cite{LaYe}, one of the
iterations of the IFS must be strictly contractive and this is
removed in this paper. It is an important improvement. Therefore,
we provide a Ruelle operator theorem for a system to which each
branch contains an indifferent fixed point (see
Remark~\ref{remark} and Example~\ref{example} in the end of this
paper).

In practice, it is difficult to calculate the spectral radius
$\varrho$ of $T$. But since $T$ is a positive operator, we have
that $\|T^{n}1\|=\| T^{n} \|$  and
$$
\varrho = \lim\limits_{n} \| T^{n} \|^{
\frac{1}{n}}=\lim\limits_{n} \| T^{n}1 \|^{ \frac{1}{n}}.
$$
Therefore, from the formula of $T^{n}1$ (see the formula before
Proposition~\ref{2.1} in \S2), a simple but useful lower bound of
$\varrho$ is
\begin{equation}~\label {15}
\min_{x \in K} \sum _{j=1}^{m} p_j(x) \leq \varrho.
\end{equation}
If we replace the $\varrho$ by $\min_{x\in K} \sum_{j=1}^m p_j(x)$
in the above theorem, we can have a simple checkable sufficient
condition.

\bigskip
\begin{cor}~\label{cor} {\em Let $(X, \{ w_{j} \}_{j=1}^{m}, \{ p_{j}
\}_{j=1}^{m})$ be a nonexpansive system. If
$$
\sup_{x \in K} \sum _{j=1}^{m} p_j (x) \cdot \sup_{y \not = x}
\frac{|w_j(x)- w_j(y)|}{|x-y|} < \min_{x \in K}\sum_{j=1}^{m}
p_{j} (x),
$$
then the Ruelle operator theorem holds for this nonexpansive
system.}
\end{cor}

\bigskip

It is obvious that if $\{ w_{j} \}_{j=1}^{m}$ is a contractive
IFS, then the conditions in the above theorem and the above
corollary and Theorem~\ref{3.5} latter are trivially satisfied.
The condition of the above theorem is similar to the {\it average
contractive} condition of Barnsley {\it et al}~\cite{BDEG} where
they assumed that $\sum_{j=1}^m p_j(x) =1$, hence $\varrho =1$. It
is also similar to the one given by Hennion~\cite{Hen}, but he
considered the case that each $p_j$ is a Lipschitz continuous
function on $X$. Regarding $T$ as defined on the Lipschitz
continuous space, he showed that the essential spectral radius
$\varrho_e(T)$ is strictly less than the spectral radius
$\varrho(T)$, and then the Ruelle operator theorem holds.
Furthermore, a general formula for the essential spectral radius
$\varrho_e (T)$ for a general $C^{\alpha}$ IFS or Zygmund IFS can
be found in~\cite{BJL}. Using this formula, one can check whether
the essential spectral radius $\varrho_e(T)$ is strictly less than
the spectral radius $\varrho(T)$, and then check the Ruelle
operator theorem. However, these methods do not work for the
weakly contractive (or, more generally, nonexpansive) case. The
reason is that, in this case, $\varrho(T)$ is not an isolated
point of the spectrum, and $\varrho
(T)=\varrho_e(T)$~(refer to~\cite{PoMa,Rug}). Note that \cite{PaPo,JiMa} contain some results showing that $\varrho
(T)=\varrho_e(T)$ is held under some weaker smoothness assumptions (for example, Dini continuity) even in the contractive case.
Therefore, the
result in this paper provides a new method to check the Ruelle
operator theorem for some weakly contractive (or, more generally,
nonexpansive) IFS.

We would like to note that most people study an IFS on some
Euclidean space. This is because the existence of a compact
invariant subset $K$ for a contractive or a weakly contractive IFS
needs the structure of a Euclidean space (see~\cite{Hut,Hat}).
However, arguments in the proofs of this paper only need to assume
that $K$ is a compact Hausdorff metric space, in particular, when
we studies a dynamical system $\sigma: K\to K$ defined on a
compact Hausdorff metric space $K$ satisfying certain Markov
property. More precisely, $K =\cup_{j=1}^{m} K_{j}$ is the union
of finitely many pairwise disjoint compact subsets
$\{K_{j}\}_{j=1}^{m}$ such that each $\sigma: K_{j}\to K$ is a
homeomorphism. Then let $w_{j}$ be the inverse of $\sigma:
K_{j}\to K$ for each $1\leq j\leq m$ and define $(K, \{
w_{j}\}_{j=1}^{m})$. It can be thought as an IFS as well. Our
results in this paper are true for such a nonexpansive IFS $(K, \{
w_{j}\}_{j=1}^{m}, \{p_{j}\}_{j=1}^{m})$.

\bigskip

The paper is organized as follows. In \S2, we will present some
elementary facts about the Ruelle operator and prove Proposition
2.1. We will introduce the Ruelle operator theorem in \S3 and set
up the basic criteria for the assertion of the Ruelle operator theorem. We will prove our main
result in \S4.

\bigskip

\section{\bf Preliminaries }

Consider the system
$$
(X, \{ w_{j} \}_{j=1}^{m}, \{ p_{j} \}_{j=1}^{m}),
$$
where $X \subseteq {\Bbb R}^{d}$ is a compact subset,  $w_j : \ X
\to X$, $1\leq j\leq m$, are continuous maps and the $p_j (x)$,
$1\leq j\leq m$, are positive functions on $X$ (they are called
weights or potentials associated with $w_j$). We say that a map $w : X\to X$ is {\it nonexpansive}
if
$$
|w (x) - w(y)| \leq |x-y|, \quad \forall x, y \in X;
$$
{\it weakly contractive} if
$$
\alpha_{w}(t) :=\sup\limits_{|x-y | \leq t}|w(x)- w(y)| < t, \quad
\forall t > 0.
$$
It is clear that contractivity implies weak contractivity which
also implies nonexpansiveness. A simple nontrivial example of a
weakly contractive map is $w(x) = x/(1+x)$ on $[0,1]$. We call
$$
(X, \{ w_{j} \}_{j=1}^{m})
$$
a weakly contractive IFS if all $w_j$, $1\leq j \leq m$, are
weakly contractive; a nonexpansive IFS if all $w_j$, $1\leq j\leq
m$, are nonexpansive.

A function $p(x)$ defined on $X$ is called Dini continuous if
$$
\int_{0}^{1} \frac{\alpha_{p}(t)}{t}  dt < \infty
$$
where
$$
\alpha_{p}(t) = \sup\limits_{|x-y | \leq t}|p(x)- p(y)|.
$$
For any $0<\theta <1$, we consider the following summation
$$
S_{\theta,p} =\sum_{n=0}^{\infty} \alpha_{p}(\theta^{n}a)
$$
where $a$ is the diameter of $X$. Then, $p(x)$ is Dini continuous
is equivalent to saying that $S_{\theta, p}$ is summable, that is,
$$
S_{\theta, p} <\infty.
$$

Throughout the paper, we always assume the potentials $p_j$'s are positive Dini continuous functions on $X$.
If $\{w_{j}\}_{j=1}^{m}$ is a
contractive IFS with the contractive constant $0<\tau<1$, that is,
$$
\sup\limits_{x\neq y\in X} \frac{|w_{j}(x)-w_{j}(y)|}{|x-y|} \leq
\tau,
$$
then the Dini condition on all $p_{j}$ can be replaced by the
summable condition
$$
\max_{1\leq j\leq m} S_{\tau,p_{j}} <\infty.
$$
However, if $\{w_{j}\}_{j=1}^{m}$ is a nonexpansive IFS, we will
not have such a constant $0<\theta<1$. Thus the Dini condition on
potentials is different from the summable condition on potentials. The methods presented before (see e.g. ~\cite
{BJL,Fan,FaJi2,FaLa,Hen,LaYe,LSV,MaUr,PrSl,Wa1,Wa2,Ye1}) do not work for the system considered in this paper. We need to find a more sharp method to prove the Ruelle operator theorem under our sufficient conditions.

\begin{de}
Let $p_j, 1 \leq j \leq m,$ be positive Dini continuous functions on $X$.
We call
$$
(X, \{ w_{j} \}_{j=1}^{m}, \{ p_{j} \}_{j=1}^{m}),
$$
a nonexpansive (or weakly contractive) system, if the IFS $(X, \{ w_j \}_{j=1}^m)$ is nonexpansive (or weakly contractive).
\end{de}

Hata studied the invariant sets of the weakly contractive IFS on
$X\subseteq {\mathbb R}^d$ in~\cite{Hat}. By using the existence
of fixed points for the weakly contractive maps, he showed the
existence of a unique nonempty compact $K \subseteq X$ invariant
under $\{w_j\}_{j=1}^{m}$, i.e.
$$
K = \bigcup _{j=1}^m w_j(K).
$$
For $J = (j_1 j_2 \cdots j_n)$,  $1 \leq j_i \leq m$,  let
$$
w_{J}(x) = w_{j_1} \circ w_{j_2} \circ \cdots \circ w_{j_n}(x).
$$
Then
$$\lim_{|J| \rightarrow \infty}|w_{J}(K)|=0
$$
and
$$
K=\bigcap_{n=1}^{\infty}\bigcup_{|J|=n}w_J(K).
$$
However, for a general IFS, an invariant set may not be unique.
However, we have

\bigskip

\begin {pro}~\label{-1}
Suppose $\{w_j\}_{j=1}^m$ is a nonexpansive IFS on the compact
subset $X$ with at least one $w_{j}$ being weakly contractive.
Then there exists a unique smallest nonempty compact set $K$ such
that
$$
K = \bigcup_{j=1}^m w_j(K).
$$
Moreover for any $x \in K$, the closure of $\{w_J(x): \ |J|= n,\;
n \in {\Bbb N} \}$ is $K$, i.e.
$$ \overline{\{w_J(x): \ |J|= n,\;
n \in {\Bbb N} \}}= K.
$$
\end{pro}

\bigskip

\noindent {\bf Proof.}  Let
$$
{\cal F} = \{ F \;|\; \ \bigcup _{j=1}^m w_j(F) \subseteq F \}.
$$
By using the standard Zorn's lemma argument, there exists a
minimal compact subset $K$ such that
$$
K = \bigcup_{j=1}^m w_j(K).
$$
To show that such $K$ is unique, we assume without loss of
generality that $w_1$ is weakly contractive. If $J_n = (1 \cdots
1)$ (n-times), then $\lim_{n \to \infty} |w_{J_n}(X)| =0$. Let
$K'$ be another minimal compact invariant set and let $x \in K$
and $y \in K'$. Then
$$
\lim_{n \to \infty} w_{J_n} (x) = \lim_{n \to \infty} w_{J_n}(y) \in K \bigcap
K'.
$$
Hence
$$
K \bigcap K' \not = \emptyset,
$$
and $w_j(K \bigcap K') \subseteq K \bigcap K'$. From the
minimality of $K$, we conclude that $K=K'$, and deduce the last
statement of the proposition. \qed

\bigskip

Throughout the paper we will consider either weakly contractive
IFS or the IFS in Proposition~\ref{-1}. Hence the set $K$ is
uniquely defined. Furthermore, we can assume without loss of
generality that the diameter
$$
|K| = \sup \{ |x-y|: x, y  \in K \} =1.
$$
Let $C(K)$ be the space of all continuous functions on $K$. For
such an system, we define an operator $T : C(K) \to C(K)$ by
$$
Tf(x) = \sum_{j=1}^{m} p_j(x) f(w_j (x)).
$$
We call $T$ the {\it Ruelle operator} assocaited to the
nonexpansive system
$$
(K, \{ w_{j} \}_{j=1}^{m}, \{ p_{j} \}_{j=1}^{m}).
$$
The dual operator $T^{\ast}$ on the measure space $M(K)$ is given
by
$$T^{\ast} \mu(E)=\sum_{j=1}^m \int_{w_j^{-1}(E)} p_j(x) d \mu(x) \quad \mbox{ for any Borel set } E \subseteq K$$
(see e.g. ~\cite{BDEG}).

For $J = (j_1 j_2 \cdots j_n)$,  $1 \leq j_i \leq m$, define
$$
w_{J}=w_{j_{1}} \circ w_{j_{2}} \circ \cdots \circ w_{j_{n}}
$$
and
$$
p_{w_{J}}(x) = p_{j_{1}}(w_{j_{2}} \circ w_{j_{3}} \circ
\cdots \circ w_{j_{n}}(x)) \cdots p_{j_{n-1}}(w_{j_{n}} (x))
p_{j_{n}}(x).
$$
Then
$$
T^{n}f(x)=\sum\limits_{|J|=n}p_{w_{J}}(x)f(w_{J}x).
$$
Let $\varrho = \varrho(T)$ be the spectral radius of $T$. Since
$T$ is a positive operator, we have that $\|T^{n}1\|=\| T^{n} \|$
and
$$
\varrho = \lim\limits_{n} \| T^{n} \|^{ \frac{1}{n}}=\lim\limits_{n} \| T^{n}1
\|^{ \frac{1}{n}}.
$$
\medskip

\begin{pro}~\label{2.1} Let $(X, \{ w_{j} \}_{j=1}^{m}, \{ p_{j} \}_{j=1}^{m})$
be a nonexpansive system with at least one weakly contractive $w_j$.
Let $T$ be the Ruelle operator on $C(K)$. Then
\begin{itemize}
\item[(i)] $\min_{x \in K} \varrho^{-n}T^{n}1(x) \leq 1 \leq
\max_{x \in K} \varrho^{-n}T^{n}1(x)$ for all  $n > 0$;
\item[(ii)] if there exist $\lambda > 0$ and $0 < h \in C(K)$ such that $T
h=\lambda h$, then $\lambda = \varrho$ and there exist $ A, B > 0$
such that
$$
A \leq \varrho^{-n}T^{n}1(x) \leq B \qquad \forall \  n > 0.
$$
\end{itemize}
\end{pro}

\medskip

\noindent{\bf Proof. } We will prove the second inequality of (i),
the first inequality  is similar. Suppose it is not true, then
there exists an integer $k$ such that
 $\|T^k 1\| < \varrho^k$. Hence
$$
\varrho  =\big(\varrho(T^k)\big)^{\frac{1}k} \leq \| T^k \|^{\frac{1}{k}} = \|T^k
1\|^{\frac 1k} <\varrho,
$$
which is a contradiction.  To prove the second assertion we let
$a_{1} = \min_{x \in K} h(x)$, $a_{2}=\max_{x \in K}h(x).$  Then
$$
0< \frac {a_{1}} {a_{2}} \leq \frac{h(x)}{a_{2}}= \frac{\lambda^{-n}}
{a_{2}}T^{n}h(x) \leq \lambda^{-n}T^{n}1(x) = \lambda^{-n} \|T^{n}\|.
$$
Similarly we can show that $\lambda^{-n} \| T^{n} \| \leq {a_{2}} /
{a_{1}}$. Hence $\varrho =\lim_{n \rightarrow
\infty}\|T^{n}\|^{\frac{1}{n}}=\lambda$. \qed

\bigskip

We call the operator $T: C(K) \to C(K)$  {\it irreducible}
(see~\cite{LaYe}) if for any non-trivial, non-negative $f \in
C(K)$ and for any $x \in K$,  there exists an integer $n>0$ such
that $T^{n}f(x)>0$.

\medskip

\begin{pro}~\label{2.2}
Let $(X, \{ w_{j} \}_{j=1}^{m}, \{ p_{j} \}_{j=1}^{m})$ be a
nonexpansive system with at least one  weakly contractive $w_j$. Then
the Ruelle operator $T$ is irreducible and
$$
{\dim} \{ h \in C(K): Th = \varrho h, \ h \geq 0 \} \leq 1.
$$
If $h\geq 0$ is a $\varrho$-eigenfunction of $T$, then $h >0$.
\end{pro}

\medskip

\noindent{\bf Proof.} The proof can be found in \cite{LaYe}. We
include the details here for the sake of completeness. For any
given $f \in C(K)$ with $f \geq 0$ and $f\not\equiv 0$, let $V =
\{ x \in K: f(x) > 0 \}$.  For any $x \in K$, by
Proposition~\ref{2.1}, there exists a multi-index $J_0$ such that
$w_{J_{0}} (x) \in V$. Let $n_{0} = |J_{0}|$, then
$$
T^{n_{0}}f(x) = \sum\limits_{|J|=n_{0}}p_{w_{J}}(x)f(w_{J}x) \geq
p_{w_{J_{0}}}(x)f(w_{J_{0}}x) > 0.
$$
This proves that $T$ is irreducible.

For the dimension of the eigensubspace, we suppose that there exist two
independent  strictly positive $\varrho$-eigenfunctions $h_1, h_2 \in C(K)$.
Without loss of generality we assume that $0< h_1 \leq h_2$ and
$h_1(x_{0})=h_2(x_{0})$
for
some $x_{0} \in K$. Then $h = h_2-h_1 (\geq 0)$ is a $\varrho$-eigenfunction of
$T$ and  $h(x_{0}) = 0$. It follows that $T^{n}h(x_{0})=\varrho^{n} h(x_{0}) = 0,$
 which contradicts to the irreducibility of $T$. Hence
the dimension of the $\varrho$-eigensubspace is at most 1.

The strict positivity of $h$ follows directly from the
irreducibility of $T$. \qed

\section{\bf Ruelle Operator Theorem}

\begin{pro}
Let $\varrho_e$ be the essential spectral radius of $T$. Suppose
$\varrho_e < \varrho$. Then  there exists a $h \in C(K)$ with
$h>0$, a probability measure $\mu \in M(K)$ and a constant $0< b
<1$ such that for any $f \in C(K)$,
$$\big\| \varrho^{-n} T^n f - \langle \mu, f \rangle h \big\|_\infty=O(b^n).$$
\end{pro}

\bigskip

\noindent {\bf Proof.} Without loss of generality, we assume that
$$\max_{x \in K} \sum_{j=1}^m p_j(x) \leq 1.$$
Then, we can prove, by induction, that
$$
\sup_{n > 0} \| T^n 1 \| =\sup_{n>0}\max_{x \in K} \sum_{|J|=n} p_{w_J}(x) \leq 1.
$$
Then, the operators sequence $n^{-1} T^n$ converges weakly to
$0$. Note that (see~\cite{Nus} or~\cite{BJL}))
$$\varrho_e=\lim_{n \to \infty} \big(\inf\{ \| T^n - Q\|: Q \mbox{ is compact on } C(K) \}\big)^{\frac{1}{n}}.$$
From this, together with the assumption $\varrho_e < \varrho$ and
theorem VIII.8.7 in~\cite{DuSc}, it follows that $T$ is
quasi-compact~\cite{Hen}. By making use of Hennion's
method~\cite{Hen}, we can deduce the assertion. \qed

\bigskip

In the following, we are interested in the case that $\varrho_e =
\varrho$. We first give a basic criterion for the existence of the
eigenfunction corresponding to the spectral radius $\varrho$ in
this case.

\medskip

\begin{pro}~\label{2.4} Let $(X, \{ w_{j} \}_{j=1}^{m}, \{ p_{j}
\}_{j=1}^{m})$ be a nonexpansive system with at least one weakly
contractive $w_j$. Suppose
\begin{itemize}
\item[(i)] there exist $A, B > 0$ such that  $A \leq \varrho^{-n}
T^{n}1(x) \leq B$ for any  $x \in K$ and $n>0$, and
\item[(ii)] for any $f \in C(K)$, $\{ \varrho^{-n} T^{n}f \}_{n=1}^{\infty}$ is an
 equicontinuous sequence.
\end{itemize}
Then there exists a unique positive function $h \in C(K)$ and a unique
probability measure $\mu \in M(K)$ such that
$$
Th = \varrho h,\ \ \ \ \ T^{*}\mu=\varrho \mu, \ \ \ \ \ \langle \mu, h \rangle=1.
$$
Moreover, for every $f \in C(K)$, $\varrho^{-n} T^{n}f$ converges to
$\langle \mu, f \rangle h$ in the supremum norm, and for every $\xi \in M(K)$,
$\varrho^{-n} T^{*n}
\xi$ converges weakly to $\langle \xi, h \rangle \mu$.
\end{pro}

\medskip

\noindent{\bf Proof.} The proof can be found in~\cite{LaYe}, and we
omit it. \qed

\medskip

\begin{de}
Let $(X,
\{w_j\}_{j=1}^m, \{p_j\}_{j=1}^m)$ be a nonexpansive system. We say that the
Ruelle operator theorem holds for this system if there exists a
unique positive function $h \in C(K)$ and a unique probability $\mu \in
M(K)$ such that
$$
Th = \varrho h,\ \ \ \ \ T^{*}\mu=\varrho \mu, \ \ \ \ \ \langle \mu, h \rangle=1,
$$
and for every $f \in C(K)$, $\varrho^{-n} T^{n}f$ converges to $\langle
\mu, f \rangle h$ in the supremum norm.
\end{de}

\medskip

In the next section, we will study the Ruelle operator theorem for
a nonexpansive system under the framework in Proposition~\ref{2.4}.

\section{\bf Some sufficient conditions }

Throughout this section we consider a nonexpansive system $
(X, \{ w_{j} \}_{j=1}^{m}, \{ p_{j} \}_{j=1}^{m})$. And, we assume the nonexpansive IFS $
(X, \{ w_{j} \}_{j=1}^{m})$ containing at least one weakly contractive $w_j$. We will prove the Ruelle operator theorem by applying
Proposition~\ref{2.4}.

In the next lemma we will see that the Dini condition on all $p_j$
also implies a similar nature property of the ``bounded distortion
property". Recall that  an equivalent condition for a function
$p(x)$ on $K$ to be Dini continuous is
$$
\sum_{n=0}^\infty \alpha_p
(\theta^n) < \infty
$$
for any $0 < \theta < 1$.

\bigskip

\begin {lemma}~\label {0}  Suppose $(X, \{ w_{j} \}_{j=1}^{m}, \{ p_{j} \}_{j=1}^{m})$
is a nonexpansive system. Let
$$
\alpha (t) = \max_{1\leq j\leq m} \{\alpha_{\log p_j}(t)\}.
$$
Let $0 < \theta <1$ and let
$$
a = \sum_{n=0}^\infty \alpha( \theta^n).
$$
For any fixed $x, y
\in K$, if $J= (j_1 \cdots j_n)$ satisfies the condition:
$$
|w_{j_{i+1}} \circ \cdots \circ w_{j_n}(x)- w_{j_{i+1}}
\circ \cdots \circ w_{j_n}(y)| \leq \theta^{n-i} \quad \forall 1 \leq i < n.
$$
Then
$$
p_{w_J}(x) \leq e^a p_{w_J}(y).
$$
\end{lemma}

\medskip

\noindent {\bf Proof.}  The inequality follows from the estimate
that
\begin{eqnarray*}
\Big| \log \frac {p_{w_J}(x)}{ p_{w_J}(y)} \Big| &
 \leq & \sum_{i=1}^n |\log p_{j_i}(w_{j_{i+1}} \circ \cdots \circ w_{j_{n}}(x)) - \log p_{j_i}(w_{j_{i+1}} \circ
\cdots \circ w_{j_{n}} (y))|\\
& \leq & \sum_{i=1}^n \alpha (\theta ^{n-i}) \leq a.
\end{eqnarray*}
\qed

\begin{pro} ~\label{3.1} Let
$(X, \{ w_{j} \}_{j=1}^{m}, \{ p_{j} \}_{j=1}^{m})$ be a
nonexpansive system. Suppose
\begin{itemize}
\item[(i)]  $r:= \sup_{x \in K} \min_{1 \leq j \leq m}\sup_{y \not=
x}\frac{|w_{j}(x)-w_{j}(y)|}{|x-y|} < 1;$
\item[(ii)] there exist constants $A, B >0$ such that $A \leq \varrho^{-n}
T^{n}1(x) \leq B$ for any  $x \in K$ and $n >0$.
\end{itemize}
Then the Ruelle operator theorem holds for this IFS.
\end{pro}

\medskip

We would like to point out that the condition (i) of Proposition \ref{3.1} is a generalization of the condition (i) of [Theorem 4.2, 15]. We extend theorem 4.2 of \cite{LaYe} so that the system considered in this paper satisfies the condition (i) of Proposition \ref{3.1}.

\bigskip

\noindent{\bf Proof.} The proof is the same as the one of [Theorem 4.2, 15], and we omit it.  \qed

\bigskip

For any integer $n$, we let $I^n=\{ J=(j_1 j_2 \cdots j_n): 1 \leq j_i \leq m \},$
and let $$D_n= \big\{ (n_1, n_2, \cdots, n_k):  0< n_i < n_{i+1} \
\mbox{ and } n_k \leq n \big\} \bigcup \{ (0) \}.$$

For any $J \in I^n$ and any $0 \leq k < l \leq n$,
we define $J|_l^k =(j_{n-l+1} j_{n-l+2} \cdots j_{n-k})$.
We let $J|_l^k =\emptyset$ if $k=l$.

For any multi-index $J$ and $x \in K$, we let
$$\gamma_J(x)=\sup_{y \not= x} \frac{|w_J(x)-w_J(y)|}{|x-y|}.$$
For convenience, we let $\gamma_J(x)=1$ and $p_{w_J}(x)=1$ if $|J|=0.$

\bigskip

\begin{pro}~\label{N-0}
Let $\{ D(k) \}_{k=1}^\ell$ be a partition of $I^n$, and let
\begin{equation} \label{N-7}
0=n_0^{(k)} < n_1^{(k)} < \cdots < n_{t_k}^{(k)}=n \ \ \ \ \forall \ 1 \leq k \leq \ell.
\end{equation}
Then for any $x \in K$,
$$\sum_{k=1}^\ell \sum_{J \in D(k)} p_{w_J}(x) \cdot \prod_{t=1}^{t_k}
\gamma_{\displaystyle{J}
\big|_{n_t^{(k)}}^{n_{t-1}^{(k)}}}\big(w_{\displaystyle{J}|_{n_{t-1}^{(k)}}^0} x\big) \leq a^n,$$
provided that
\begin{equation}~\label{N-1}
\sup_{x \in K} \sum_{j=1}^m p_j (x) \cdot \gamma_j(x) \leq a.
\end{equation}
\end{pro}

\bigskip

\noindent {\bf Proof.} Note the fact that
$$p_{w_J}(x)= \prod_{i=0}^{n-1} p_{j_{n-i}} \big(w_{\displaystyle{J}|_i^0} x\big) \quad \forall \ J=(j_1 j_2
\cdots j_n).$$ From (\ref{N-1}), we can deduce inductively that
for any integer $n$,
\begin{equation}~\label{N-8}
\sum_{|J|=n} p_{w_J}(x) \cdot \prod_{i=0}^{n-1} \gamma_{J|_{i+1}^i} (w_{J|_i^0}x) \leq a^n.
\end{equation}

For any multi-index $J=(j_1 j_2
\cdots j_N)$ and $x \in K$, we have
$$\frac{|w_J(x)-w_J(y)|}{|x-y|} = \prod_{i=0}^{N-1}
\frac{|w_{j_{N-i}}(w_{J|_{i}^0}x)-w_{j_{N-i}}(w_{J|_{i}^0}y)|}{|w_{J|_{i}^0}(x)-w_{J|_{i}^0}(y)|},
\quad \forall y \not= x.$$
This implies that
\begin{equation}~\label{N-2}
\gamma_J (x) \leq \prod_{i=0}^{N-1} \gamma_{j_{N-i}} (w_{J|_i^0}x).
\end{equation}
From the assumption (\ref{N-7}), using the same argument as
(\ref{N-2}), we deduce that for any $J$ with $|J|=n$,
\begin{equation} \label{N-9}
\prod_{t=1}^{t_k} \gamma_{\displaystyle{J} \big|_{n_t^{(k)}}^{n_{t-1}^{(k)}}}
\big(w_{\displaystyle{J}|_{n_{t-1}^{(k)}}^0} x\big) \leq \prod_{i=0}^{n-1} \gamma_{J|_{i+1}^i} (w_{J|_i^0}x).
\end{equation}

Note that $\{ D(k) \}_{k=1}^\ell$ is a partition of $I^n ( =\{ J: |J|=n \})$. We have
\begin{eqnarray*}
&& \sum_{k=1}^\ell \sum_{J \in D(k)} p_{w_J}(x) \cdot \prod_{t=1}^{t_k}
\gamma_{\displaystyle{J} \big|_{n_t^{(k)}}^{n_{t-1}^{(k)}}}
\big(w_{\displaystyle{J}|_{n_{t-1}^{(k)}}^0} x\big)\\
& \leq & \sum_{|J|=n} p_{w_J}(x) \cdot \prod_{i=0}^{n-1}
\gamma_{J|_{i+1}^i} (w_{J|_i^0}x) \qquad \mbox{(by (\ref{N-9}))}\\
& \leq & a^n \qquad \mbox{(by (\ref{N-8}))}.
\end{eqnarray*}
Thus, the conclusion follows.
\qed

\bigskip

As a consequence of Proposition~\ref{3.1}, we have

\bigskip

\begin{pro}~\label{3.5-1} Let $(X, \{ w_{j} \}_{j=1}^{m}, \{ p_{j}
\}_{j=1}^{m})$ be a nonexpansive system. Suppose that
\begin{itemize}
\item[(i)] there
exists $k$ such that
$$\sup_{x \in K} \sum_{|J|=k} p_{w_{J}}(x) \cdot \gamma_J (x) <
\varrho^{k};$$
\item[(ii)] there exist constants $A, B >0$ such that $A \leq \varrho^{-n}
T^{n}1(x) \leq B$ for any  $x \in K$ and $n >0$.
\end{itemize}
Then the Ruelle operator theorem holds.
\end{pro}

\bigskip

\noindent{\bf Proof.} By (i) there exists a $0< \eta <1$ such that
$$\sup_{x \in K} \sum_{|J|=k} p_{w_{J}}(x) \cdot \gamma_J (x) \leq \eta \varrho^{k}$$
This, together with Proposition \ref{N-0}, implies that for any $x \in K$ and $\ell \in {\Bbb{N}}$,
$$\sum_{|J|= \ell k} p_{w_J}(x) \cdot \prod_{t=1}^{\ell}
\gamma_{\displaystyle{J}
\big|_{t k}^{(t-1) k}} \big(w_{\displaystyle{J}|_{(t-1) k}^0} x\big) \leq \eta^\ell \varrho^{\ell k}.$$
By using the argument similar to (\ref{N-2}), we can prove that for any muti-index $J$ with $|J|=\ell k$,
$$\gamma_J(x) \leq \prod_{t=1}^{\ell}
\gamma_{\displaystyle{J}
\big|_{t k}^{(t-1) k}} \big(w_{\displaystyle{J}|_{(t-1) k}^0} x\big).$$
It follows that
\begin{equation}~\label{N-3}
 \sum_{|J|= \ell k} p_{w_J}(x) \cdot \gamma_J(x) \leq \eta^\ell \varrho^{\ell k}.
\end{equation}
We claim that
$$\sup_{x \in K}\inf_{\ell \in {\Bbb{N}}}\min_{|J|=\ell k}\gamma_J(x) =0.$$
Otherwise, we suppose that
$$\sup_{x \in K}\inf_{\ell \in {\Bbb{N}}}\min_{|J|=\ell k}\gamma_J(x) >0.$$
Then, there exists a $b_0 >0$ and a $x_0 \in K$ such that
$$\inf_{\ell \in {\Bbb{N}}}\min_{|J|=\ell k}\gamma_J(x_0) \geq b_0.$$
This, combined with (\ref{N-3}) and (ii), implies that for any $\ell \in {\Bbb{N}}$,
\begin{eqnarray*}
\eta^\ell & \geq & \varrho^{- \ell k} \sum_{|J|= \ell k} p_{w_J}(x_0) \cdot \gamma_J(x_0) \geq b_0 \cdot \varrho^{- \ell k} \sum_{|J|= \ell k} p_{w_J}(x_0)\\
& = & b_0 \cdot \varrho^{- \ell k} T^{\ell k}1(x_0) \geq b_0 A. \qquad \mbox{ (by (ii))}
\end{eqnarray*}
This contradicts to the choice of $0< \eta <1$. Then, the claim follows. And thus, there exists a $\ell_0 \in {\Bbb{N}}$ and a $J_0$ with $|J_0|=\ell_0 k$ such that $\sup_{x \in K}\gamma_{J_0}(x)<1.$
Hence, by Proposition \ref{3.1}, the Ruelle operator theorem for $T^{\ell_0 k}$ holds. This implies that the Ruelle operator theorem for $T$ holds. \qed

\bigskip

\begin{thm}~\label{3.5} Suppose $(X, \{ w_{j} \}_{j=1}^{m}, \{ p_{j}
\}_{j=1}^{m})$ is a nonexpansive system. If there
exists $k$ such that
\begin{equation}~\label{31}
\sup_{x \in K} \sum_{|J|=k} p_{w_{J}}(x) \cdot \gamma_J (x) <
\varrho^{k},
\end{equation}
then the Ruelle operator theorem holds.
\end{thm}

\bigskip

\noindent{\bf Proof.} Since the Ruelle operator theorem for
$T^{k}$ implies the Ruelle operator theorem for $T$, we may assume
$k=1$ in the hypothesis, so that (\ref{31}) is reduced to
\begin{equation}~\label{20}
\sup_{x \in K} \sum_{j=1}^{m}p_{j}(x) \cdot \gamma_j (x) <
\varrho.
\end{equation}
This means that the condition (i) of of Proposition~\ref{3.5-1} is satisfied. Hence, we need only to show that condition
(ii) of Proposition~\ref{3.5-1} is also satisfied, i.e. there exist $A, B
>0$ such that
$$
A \leq \varrho^{-n}\sum\limits_{|J|=n}p_{w_{J}}(x) \leq B \quad \forall \ n.
$$
By (\ref{20}) we can find $0<\eta<1$ such that
\begin{equation} \label{N-6}
\sup_{x \in K} \sum_{j=1}^{m}p_{j}(x) \cdot \gamma_j (x) \leq
\eta\varrho.
\end{equation}

For any fixed $x \in K$, choose $\theta$ such that $0<\eta<\theta <
1$. For any integer $n$ and $J \in I^n$, let $n_1$ be the largest integer such that
$$\gamma_{\displaystyle{J}|_{n_1}^0}(x) \geq \theta^{n_1},$$
and let $n_2 ( > n_1)$ be the largest integer such that
$$\gamma_{\displaystyle{J}|_{n_2}^{n_1}}(w_{\displaystyle{J}|_{n_1}^0}x) \geq \theta^{n_2-n_1},$$
and so on. Then, we find a sequence $\{ n_i \}_{i=1}^{t_J}$ such that
$$\gamma_{\displaystyle{J}|_{n_{i+1}}^{n_i}}(w_{\displaystyle{J}|_{n_i}^0}x)
\geq \theta^{n_{i+1}-n_i} \quad \forall  \ 1 \leq i \leq n_{t_J} -1,$$
and
\begin{equation} ~\label{N-10}
\gamma_{\displaystyle{J}|_i^{n_{t_J}}}(w_{\displaystyle{J}|_{n_{t_J}}^0}(x))
< \theta^{i - n_{t_J}} \quad \forall \ n_{t_J} < i \leq n.
\end{equation}
Define $\sigma: I^n \to D_n$ by
$$\sigma(J)=(n_1, n_2, \cdots, n_{t_J}).$$
Then $\# \sigma(I^n) < \infty.$ Denote $\sigma(I^n)=\{ A_k \}_{k=1}^\ell,$ where $A_k \in D_n$. Let
$$D(k)=\{ J: \sigma(J)=A_k \}, \ \ \forall 1 \leq k \leq \ell.$$
It is clear that
$$D(i) \bigcap D(j) =\emptyset, \quad \forall \ i \not= j.$$
Hence, $\{ D(k) \}_{k=1}^\ell$ is  a partition of $I^n$.

For any $1 \leq k \leq \ell$, let $A_k=(n_1^{(k)}, n_2^{(k)},
\cdots, n_{t_{k-1}}^{(k)})$. For convenience, we let $n_0^{(k)}=0$ and let $n_{t_k}^{(k)}=n$.
By making use of (\ref{N-6}), it follows from Proposition \ref{N-0} that
\begin{eqnarray} \label{6}
S_0: = \sum_{k=1}^\ell \sum_{J \in D(k)} p_{w_J}(x)
\cdot \prod_{t=1}^{t_k} \gamma_{\displaystyle{J} \big|_{n_t^{(k)}}^{n_{t-1}^{(k)}}}
\big(w_{\displaystyle{J}|_{n_{t-1}^{(k)}}^0} x\big) \leq (\eta \varrho)^n.
\end{eqnarray}

Let
\begin{eqnarray*}
\Omega(n, k)&=&\{ J: |J|=n \mbox{ and } n_{t_J} = k \}, \quad 1\leq k \leq n, \\
\Omega(n, 0)&=&\{ J: |J|=n \mbox{ and }  n_{t_J}=0 \}.
\end{eqnarray*}
Then
$$
I^n = \bigcup_{k=0}^{n}\Omega(n, k).
$$

Without loss of generality, we assume that $\Omega(n, n)=\big\{ D(k) \big\}_{k=1}^{\ell_0},$
where $\ell_0 \leq \ell.$ And we let
$$S_1:= \sum_{k=1}^{\ell_0} \sum_{J \in D(k)}p_{w_J}(x)
\cdot \prod_{t=1}^{t_k} \gamma_{\displaystyle{J}
\big|_{n_t^{(k)}}^{n_{t-1}^{(k)}}}\big(w_{\displaystyle{J}|_{n_{t-1}^{(k)}}^0} x\big).$$

For any $1 \leq k \leq \ell_0$ and any $J \in D(k)$, we have
$n_{t_{k-1}}^{(k)}=n_{t_J}=n$, and this implies that
$$\prod_{t=1}^{t_k} \gamma_{\displaystyle{J} \big|_{n_t^{(k)}}^{n_{t-1}^{(k)}}}
\big(w_{\displaystyle{J}|_{n_{t-1}^{(k)}}^0} x\big)
\geq \prod_{t=1}^{t_k} \theta^{n_t^{(k)} - n_{t-1}^{(k)}} = \theta^n.$$
From this, we conclude that
$$S_1 \geq \sum_{k=1}^{\ell_0} \sum_{J \in D(k)}p_{w_J}(x)
\cdot \theta^n = \sum_{J \in \Omega(n, n)} p_{w_J}(x) \cdot \theta^n.$$
This, combined with (\ref{6}), implies that
$$\sum_{J \in \Omega(n, n)} p_{w_{J}}(x) \cdot \theta^n \leq S_1 \leq S_0 \leq (\eta \varrho)^n.$$
Thus, it follows that
\begin{equation}~\label{8}
\varrho^{-n}\sum_{J \in \Omega(n, n)} p_{w_{J}}(x) \leq
(\frac{\eta}{\theta})^n.
\end{equation}

Remember that
$$
\alpha(t)=\max_{1 \leq j \leq m}\alpha_{\log
p_{j}}(t)
$$
and
$$
a=\sum_{k=0}^{\infty} \alpha(\theta^{k}).
$$
Then $a$ is finite because all the $p_{i}$ are Dini
continuous functions on $X$. For any $n>0$, we can make use of
Proposition~\ref {2.1}(i) to find $x_{n} \in K$ such that
\begin{equation}~\label{9}
\varrho^{-n}\sum_{|J|=n} p_{w_{J}}(x_{n}) \leq 1.
\end{equation}
For any $J=(j_1 j_2 \cdots j_n) \in \Omega(n, k)$, we have $J|_k^0
\in \Omega(k, k)$. By using (\ref{N-10}), we can deduce from
Lemma~\ref{0} that
$$p_{w_{J|_n^k}}(w_{J|_k^0}x) \leq e^a p_{w_{J|_n^k}}(y) \quad \forall \ y \in K.$$
(We use $|K|=1$ here.) Hence
\begin{eqnarray}~\label{10}
p_{w_{J}}(x)=p_{w_{J|_n^k}}(w_{J|_k^0}x) p_{w_{J|_k^0}}(x)
 \leq  e^{a} p_{w_{J|_n^k}}(x_{n-k}) p_{w_{J|_k^0}}(x).
\end{eqnarray}
It follows that
\begin{eqnarray}~\label{11}
&& \varrho^{-n}\sum_{|J|=n} p_{w_{J}}(x)
= \varrho^{-n}\sum_{k=0}^{n}\sum_{J \in \Omega(n, k)}p_{w_{J}}(x) \nonumber\\
&  \leq& \varrho^{-n}\sum_{k=0}^{n}\sum_{J \in \Omega(n,
k)}e^{a}p_{w_{J|_n^k}}(x_{n-k}) p_{w_{J|_k^0}}(x)  \qquad (\mbox{by }\
(\ref{10}))\nonumber\\
&  \leq& e^{a} \sum_{k=0}^{n}\Big(\varrho^{-n+k} \sum\limits_{ |J^{'}| =n-k }
p_{w_{J^{'}}}(x_{n-k})\Big) \Big(\varrho^{-k}\sum\limits_{J^{''} \in
\Omega(k, k)} p_{w_{J^{''}}}(x)\Big)\nonumber\\
&  \leq & e^{a} \sum\limits_{k=0}^{n}1 \cdot (\frac{\eta}{\theta})^k
 \qquad (\mbox{by }(\ref{8}), \ (\ref{9})).
\end{eqnarray}
The last term is bounded by $e^{a}
\sum_{k=0}^{\infty}(\frac{\eta}{\theta})^{k}:=B_{1}$. This concludes the upper
bound estimate.

 For the lower bound estimation, we note that Proposition~\ref{2.1}(i) and (\ref{11}) implies that for any
$n>0$, there exists $y_{n} \in K$ such that
$$
1 \leq C_{n}:=\varrho^{-n}\sum\limits_{|J|=n} p_{w_{J}}(y_{n}) \leq B_{1}.
$$

For any fixed $x \in K$, we let
$$\alpha_J
=\sum_{i=0}^{n-1}\alpha(|w_{J|_i^0}(x)-w_{J|_i^0}(y_n)|).$$
Then, we have
 $$p_{w_J}(y_n) \leq p_{w_J}(x) e^{\alpha_J}.$$
By (\ref{N-10}), we have for any $J \in \Omega(n, k)$,
$$|w_{J|_i^0}(x)-w_{J|_i^0}(y_n)| < \theta^{i - k} \quad \forall \ k < i \leq n.
$$
(We use $|K|=1$ here.) It follows that
$$\alpha_{J} \leq a + k \alpha(1) \quad \forall J \in \Omega(n, k).$$

Using the same argument as (\ref{11}), we can deduce that
$$\varrho^{-n} \sum_{J \in \Omega(n,k)}p_{w_{J}}(y_{n}) \leq e^a (\frac{\eta}{\theta})^k.$$
And then, we have
\begin{eqnarray*}
 &&\varrho^{-n}\sum_{|J|=n} p_{w_{J}}(y_{n}) \alpha_{J}
=\varrho^{-n} \sum_{k=0}^{n}\sum_{J \in \Omega(n,k)}p_{w_{J}}(y_{n}) \alpha_J\\
&\leq & \varrho^{-n} \sum_{k=0}^{n}\big( a+ k \alpha(1)\big)\sum_{J \in
\Omega(n,k)} p_{w_{J}}(y_{n})\\
& \leq & e^a \sum_{k=0}^{n}\big( a+ k \alpha(1)\big) (\frac{\eta}{\theta})^k \leq B_2,
\end{eqnarray*}
where $B_2:=e^a \sum_{k=0}^\infty \big( a+ k \alpha(1)\big) (\frac{\eta}{\theta})^k$.
By the convexity of function $e^{x}$, we have
\begin{eqnarray*}
\varrho^{-n}\sum_{|J|=n} p_{w_J}(x) &\geq &\varrho^{-n}\sum_{|J|=n}
p_{w_{J}}(y_{n}) e^{- \alpha_{J}} \geq \frac{\varrho^{-n}}{C_{n}}\sum_{|J|=n}
p_{w_{J}}(y_{n}) e^{- \alpha_{J}}\\
&\geq& e^{-\frac{1}{C_{n}}\varrho^{-n}\sum_{|J|=n} p_{w_{J}}(y_{n}) \alpha_{J}}
\geq  e^{-B_{2}}.
\end{eqnarray*}
This completes the proof. \qed

\bigskip
\begin{re}~\label{remark}
We note that for any muti-index $J$ and $x \in K$,
$$\gamma_J (x) \leq \sup_{y \not=
z}\frac{|w_J(z)-w_J(y)|}{|z-y|}.
$$
And then, for any integer $n$, we have
$$
\sum_{|J|=n} p_{w_{J}}(x) \cdot \gamma_J (x) \leq \sum_{|J|=n} p_{w_{J}}(x) \cdot \sup_{y \not=
z}\frac{|w_J(z)-w_J(y)|}{|z-y|}.
$$
Hence, Theorem~\ref{3.5} in this paper is a generalization of
theorem 4.4 in~\cite{LaYe}. However, the following example
indicates that this generalization is non-trivial.
\end{re}

\bigskip
\begin{ex}~\label{example}
Let $X=[0, 1]$, and let $w_1(x)= x- \frac{x^2}{2}$, $w_2(x)=\frac
{1}{2} + \frac{x^2}{2}$. Then $w_1'(\cdot) \geq 0$, $w_2'(\cdot)
\geq 0$. And $w_1(0)=0$, $w_2(1)=1$; $w_1'(0)=w_2'(1)=1$ and
$w_1'(1)=w_2'(0)=0$.
\end{ex}

In this example, both $w_1$ and $w_2$ are not strictly contractive.
In fact, $0$ is the indifferent fixed point of $w_1$; and $1$
is the indifferent fixed point of $w_2$. It is easy to see that
the IFS $(X, \{ w_j \}_{j=1}^2)$ is weakly contractive.

Let $p_1$ be any positive Dini function (not a Lipschitz function)
on $X$ with the inequalities $0<p_1 (\cdot) <1.$
Let
$$\delta =\frac{1}{5} \cdot \min_{x \in X} \{ p_1(x), 1-p_1(x) \}>0,$$
 and let
$$
g(x) = \left\{\begin{array}{cc}
\delta - 2^{-1} + x , &\ {\text {if}}\ \ 2^{-1} - \delta < x \leq 2^{-1}\\
\delta + 2^{-1} - x, & \ {\text {if} }\ \ 2^{-1} < x < 2^{-1} + \delta \\
0, & \ {\text {otherwise} }.
\end{array}
\right.
$$
Define a Dini function $p_2$ on $X$ by
$$p_2 (x)=1 - p_1 (x) + g(x) \quad \forall x \in X.$$
Then
$$1 \leq \sum_{j=1}^2 p_j (x) = 1 +g(x) \leq 1+\delta.$$
And for any $x \in X$,
\begin{equation} ~\label{N-4}
g(x)-\frac{1}{4} p_1(x) < 0 \ \ \mbox{ and } \ \ g(x) - \frac{1}{4} p_2(x) < 0.
\end{equation}
Let $K$ be the invariant set of the IFS $\{ w_j \}_{j=1}^2$. Define
$$T f(x) = \sum_{j=1}^2 p_j (x) f \circ w_j(x), \quad \forall \ f \in C(K).$$
Let $\varrho$ be the spectral radius of the operator $T$.
Then, we have
\begin{equation} ~\label{N-5}
1 \leq \varrho \leq 1+\delta.
\end{equation}
Note that
\begin{eqnarray*}
\gamma_1(x)& = & \sup_{y \not= x} \frac{|w_1(y)-w_1(x)|}{|y-x|} = \sup_{y \not= x} \frac{|y-x - 2^{-1}y^2 + 2^{-1}x^2|}{|y-x|}\\
& = & \sup_{y \not= x} \big(1- \frac 1 2 (x+y)\big)= 1 - \frac {x} {2};\\
\gamma_2(x) & = & \sup_{y \not= x} \frac{|w_2(y)-w_2(x)|}{|y-x|}=\sup_{y \not= x} \frac{|2^{-1}y^2 - 2^{-1}x^2|}{|y-x|}=\frac 1 2 (1+x).
\end{eqnarray*}
We have
\begin{eqnarray*}
\sum _{j=1}^2 p_j (x) \cdot \gamma_j (x) & = & p_1 (x) \cdot (1 - \frac {x} {2})
+ p_2 (x) \cdot \frac {1+x} {2}\\
& \leq & \left\{\begin{array}{cc}
p_1(x)+ \frac{3}{4} p_2(x), &{\text {if}}\ \ 0 \leq x \leq \frac{1}{2},\\
\frac{3}{4} p_1(x) + p_2(x), & {\text {if}}\ \ \frac{1}{2} < x \leq 1,
\end{array}
\right.\\
& \leq & \left\{\begin{array}{cc}
1 + \big(g(x) - \frac{1}{4} p_2(x)\big), &{\text {if}}\ \ 0 \leq x \leq \frac{1}{2},\\
1 + \big(g(x) - \frac{1}{4} p_1(x)\big), & {\text {if}}\ \ \frac{1}{2} < x \leq 1,
\end{array}
\right.\\
& < & 1. \qquad \mbox{ (by (\ref{N-4})) }
\end{eqnarray*}
This, together with (\ref{N-5}), implies that
$$
\sup_{x \in X} \sum _{j=1}^2 p_j (x) \cdot \gamma_j (x) < \varrho.
$$
And then, Theorem~\ref{3.5} implies that the Ruelle operator theorem holds
for this weakly contractive system.

Because of the equalities
$$\sup_{y \not= z} \frac{|w_j(y)-w_j(z)|}{|y-z|}=1 \ \ \ \forall j=1, 2,$$
and
$$\sup_{x \in X} \sum _{j=1}^2 p_j (x) = \sup_{x \in X} \big(1 + g(x)\big) = 1+\delta,$$
by noting that (\ref{N-5}), the following inequality:
$$
\sup_{x \in X} \sum _{j=1}^2 p_j (x) \cdot 1 < \varrho
$$
does not hold. Hence, for this system, the condition of theorem
1.2 in~\cite{LaYe} is not satisfied.

\bigskip
\bigskip

\noindent Yunping Jiang: Department of Mathematics, Queens College
of CUNY, Flushing, NY 11367 and Department of Mathematics, CUNY
Graduate Center, 365 Fifth Avenue, New York, NY 10016.

\noindent E-mail: Yunping.Jiang@@qc.cuny.edu

\medskip

\noindent Yuan-Ling Ye: School of Mathematical Sciences, South
China Normal University, Guangzhou 510631, People's Republic of
China

\noindent E-mail: ylye@@scnu.edu.cn

\bigskip

\end{document}